\def\E{\end{document}}
\documentclass[11pt]{article}
\usepackage{amssymb,amsmath}

\begin{document}
\title{On the Global Existence and Blowup Phenomena of Schr\"{o}dinger Equations with Multiple Nonlinearities}
 \author{Xianfa Song {\thanks{E-mail: songxianfa2004@163.com(or
 songxianfa2008@sina.com)
}}\\
\small Department of Mathematics, School of Science, Tianjin University,\\
\small Tianjin, 300072, P. R. China }

\maketitle
\date{}

\newtheorem{theorem}{Theorem}[section]
\newtheorem{definition}{Definition}[section]
\newtheorem{lemma}{Lemma}[section]
\newtheorem{proposition}{Proposition}[section]
\newtheorem{corollary}{Corollary}[section]
\newtheorem{remark}{Remark}[section]
\renewcommand{\theequation}{\thesection.\arabic{equation}}
\catcode`@=11 \@addtoreset{equation}{section} \catcode`@=12

\begin{abstract}

In this paper, we consider the global existence and blowup
phenomena of the following Cauchy problem
\begin{align*} \left\{\begin{array}{ll}&-i
u_t=\Delta u-V(x)u+f(x,|u|^2)u+(W\star|u|^2)u, \quad x\in\mathbb{R}^N, \quad t>0,\\
&u(x,0)=u_0(x), \quad x\in\mathbb{R}^N,
\end{array}
\right. \end{align*} where $V(x)$ and $W(x)$ are real-valued
potentials with $V(x)\geq 0$ and $W$ is even, $f(x,|u|^2)$
is measurable in $x$ and continuous in $|u|^2$,
and $u_0(x)$ is a complex-valued function of $x$. We
obtain some sufficient conditions and establish two sharp
thresholds for the blowup and global existence of the solution to
the problem. These results can be looked as the supplement to Chapter 6 of \cite{Cazenave2}.
In addition, our results extend those of \cite{Zhang} and improve some of \cite{Tao2}.



\end{abstract}

\section{Introduction}
\qquad In this paper, we are interested in the global existence and
blowup phenomena of the following Cauchy problem
\begin{align} \label{system0}\left\{\begin{array}{ll}&-i
u_t=\Delta u-V(x)u+f(x,|u|^2)u+(W\star|u|^2)u, \quad x\in\mathbb{R}^N, \quad t>0,\\
&u(x,0)=u_0(x)\in \Sigma, \quad x\in\mathbb{R}^N,
\end{array}
\right. \end{align} where $V(x)$ and $W(x)$ are real-valued
potentials with $V(x)\geq 0$ and $W$ is even, $f(x,|u|^2)$
is measurable in $x$ and continuous in $|u|^2$,
and $u_0(x)$ is a complex-valued function of $x$, and $\Sigma$ is a natural Hilbert space:
\begin{align}
\Sigma=\{u\in H^1(\mathbb{R}^N):
\int_{\mathbb{R}^N}V(x)|u|^2dx<+\infty\}\label{97236}
\end{align}
with the inner product
\begin{align}
<\varphi,\psi>=\int_{\mathbb{R}^N}[\varphi\bar{\psi}+\nabla
\varphi\cdot \nabla \bar{\psi}+V(x)\varphi\bar{\psi}]dx\label{97237}
\end{align}
and the norm
\begin{align}
\|u\|_{\Sigma}^2=\int_{\mathbb{R}^N}[|u|^2+|\nabla
u|^2+V(x)|u|^2]dx.\label{97238}
\end{align}
The model (\ref{system0}) appears in the theory of Bose-Einstein
condensation, nonlinear optics and theory of water waves
(see\cite{Cazenave2, Ginibre, Ginibre2, Hitoshi,Kato, Oh}).

In convenience, we will give some assumptions on $V$, $f$ and $W$ as follows.

(V1)  $V(x)\geq 0$ and $V\in L^r(\mathbb{R}^N)+L^{\infty}(\mathbb{R}^N)$ for $r\geq
1$, $r>\frac{N}{2}$ or

(V2) $V(x)\in {\bf S}^c_1$, $V(x)\geq 0$ and $|D^{\alpha}V|$
is bounded for all $|\alpha|\geq 2$. Here
${\bf S}^c_1$ is the complementary set of ${\bf S}_1=\{V(x)\ {\rm  satifies } \ (V1)\}$.

(f1) $f: \mathbb{R}^N\times\mathbb{R}\rightarrow \mathbb{R}$ is measurable in $x$
and continuous in $|u|^2$ with $f(x,0)=0$.  Assume that for every $k>0$
there exists $L(k)<+\infty$ such that $|f(x,s_1)-f(x,s_2)|\leq
L(k)|s_1-s_2|$ for all $0\leq s_1< s_2<k$. Here
\begin{align} \label{98101}\left\{\begin{array}{ll}&L(k)\in C([0,\infty)),\qquad {\rm if}\ N=1\\
&L(k)\leq C(1+k^{\alpha}) \quad {\rm with}\ 0\leq
\alpha<\frac{2}{N-2}, \qquad {\rm if}\ N\geq 2.
\end{array}
\right. \end{align}

(W1) $W$ is even and $W\in
L^q(\mathbb{R}^N)+L^{\infty}(\mathbb{R}^N)$ for some $q\geq 1$,
$q>\frac{N}{4}$.

Denote $\frac{1}{(N-2)^+}=+\infty$ when $N=1,2$
and $(N-2)^+=N-2$ when $N\geq 3$.

First, we consider the local well-posedness of (\ref{system0}).
We have a proposition as follows.

{\bf Proposition 1.1.} (Local Existence Result) {\it Assume that $(f1)$ and $(W1)$ are true, $V(x)$ satisfies $(V1)$ or $(V2)$,  $u_0\in \Sigma$. Then there exists a unique
solution $u$ of (\ref{system0}) on a maximal time interval
$[0,T_{\max})$ such that $u\in C(\Sigma;[0,T_{\max}))$ and either
$T_{\max}=+\infty$ or else
$$T_{\max}<+\infty, \quad \lim_{t\rightarrow T_{\max}}\|u(\cdot,t)\|_{\Sigma}=+\infty.$$}

 {\bf Definition 1.1.} {\it If $u\in C(\Sigma;
[0,T))$ with $T=\infty$, we say that the
solution $u$ of (\ref{system0}) exists globally. If
$u\in C(\Sigma;[0,T))$ with $T<+\infty$
and $\lim_{t\rightarrow T}\|u(\cdot,t)\|_{\Sigma}\rightarrow +\infty$, we say that the solution
$u$ of (\ref{system0}) blows up in finite time.}

Our main topic is the global existence and blowup
phenomena of the solution to (\ref{system0}), which is directly motivated by \cite{Cazenave2}. Since Cazevave established some results on blowup and global existence of the solutions to
(\ref{system0}) with (V1), (f1) and (W1) in \cite{Cazenave2}, we are interested
in the parallel problems such as: What are the results about the blowup and global existence of the solutions to
(\ref{system0}) with (V2), (f1) and (W1)? How can we establish
the sharp threshold for global existence and blowup of the solution to
(\ref{system0})?

About the topic of global existence and blowup in finite time, there are many results on the special cases of (\ref{system0}). However, we only cite some very related references which only gave some sufficient conditions on global existence and blowup of the solution to the special case of (\ref{system0}). We will show how all the cited results give coherence and connection to our paper below. A special case of (\ref{system0}) is
\begin{align} \label{972310}\left\{\begin{array}{ll}&-i
u_t=\Delta u+f(|u|^2)u, \quad x\in\mathbb{R}^N, \ t>0,\\
&u(x,0)=u_0(x), \quad x\in\mathbb{R}^N.
\end{array}
\right. \end{align}  In \cite{Glassey}, Glassey established some
blowup results for (\ref{972310}). In \cite{Berestycki}, Berestyki
and Cazenave established the sharp threshold for blowup of
(\ref{972310}) with supercritical nonlinearity by considering a
constrained variational problem. In \cite{Weinstein}, Weinstein
presented a relationship between the sharp criterion for the
global solution of (\ref{972310}) and the best constant in the
Gagliardo-Nirenberg's inequality. In \cite{Cazenave3}, Cazenave
and Weisseler established the local existence and uniqueness of
the solution to (\ref{972310}) with $f(|u|^2)u=|u|^{\frac{4}{N}}u$.
Very recently, Tao et al. in \cite{Tao2}
studied the Cauchy problem (\ref{972310}) with $f(|u|^2)u=
\mu|u|^{p_1}u+\nu|u|^{p_2}u$, where $\mu$ and $\nu$ are real
numbers, $0<p_1<p_2<\frac{4}{N-2}$ with $N\geq 3$. This type of
nonlinearity brings the failure of the equation in (\ref{972310})
to be scale invariant and it cannot satisfy the conditions of the blowup theorem in
\cite{Glassey} in some cases. Tao et al. established the results on local and
global well-posedness, asymptotic behavior (scattering) and finite
time blowup under some assumptions. These papers
above have given some sufficient conditions on global
existence and blowup of the solution or established the sharp
threshold for the special case of (\ref{972310}). Naturally, we
want to establish a new sharp threshold for global existence and
blowup of the solution to (\ref{972310}) in this paper, which will
generalize or even improve these results above.

The following Cauchy problem
\begin{align} \label{972311}\left\{\begin{array}{ll}&-i
u_t=\frac{1}{2}\Delta u-V(x)u+|u|^pu, \quad x\in\mathbb{R}^N, \quad t>0,\\
&u(x,0)=u_0(x), \quad x\in\mathbb{R}^N
\end{array}
\right. \end{align}  is also a special case of (\ref{system0}). If $p<\frac{4}{N}$,
in \cite{Oh}, Oh obtained the local well-posedness and global
existence results of (\ref{972311}) under some conditions on $V(x)$. If $\frac{4}{N}\leq p<\frac{4}{(N-2)^+}$, in \cite{Zhang},
Zhang established a sharp threshold for the global existence and
blowup of the solutions to (\ref{972311}) with $V(x)=|x|^2$.
Another special case of (\ref{system0}) is the following
Cauchy problem of Schr\"{o}dinger-Hartree equation:
\begin{align} \label{9829w1}\left\{\begin{array}{ll}&-i
u_t=\Delta u+(W\star|u|^2)u, \quad x\in\mathbb{R}^N, \quad t>0,\\
&u(x,0)=u_0(x), \quad x\in\mathbb{R}^N,
\end{array}
\right. \end{align}  Using a contraction mapping argument and
energy estimates, Hitoshi obtained
the local and global existence results on (\ref{9829w1}) in \cite{Hitoshi}.
More recently, Miao et al. studied the global well-posedness and scattering for
the mass-critical Hartree equation with radial data in
\cite{Miao1} and global well-posedness, scattering and blowup for
the energy-critical, focusing Hartree equation with the radial
case in \cite{Miao2}. And in \cite{Li}, Li et al. also dealt with
the focusing energy-critical Hartree equation, they prove that the
maximal-lifespan $I=\mathbb{R}$, moreover, the solution scatters
in both time directions. However, there are few results on the sharp
threshold for global existence and blowup of
the solution to (\ref{9829w1}). Therefore, we want to establish
a sharp threshold for global existence and blowup of the solution to (\ref{9829w1})
under some conditions.

Now we will introduce some notations. Denote
\begin{align}
&F(x,|u|^2)=\int_0^{|u|^2}f(x,s)ds,\quad
G(|u|^2)=\frac{1}{4}\int_{\mathbb{R}^N}
(W\star|u|^2)|u|^2dx\label{97255}\\
&h(u)=-V(x)u+f(x,|u|^2)u+(W\star|u|^2)u,\label{98121}\\
&H(u)=-\frac{1}{2}\int_{\mathbb{R}^N}
V(x)|u|^2dx+\frac{1}{2}\int_{\mathbb{R}^N}F(x,|u|^2)dx+\frac{1}{4}\int_{\mathbb{R}^N}
(W\star|u|^2)|u|^2dx.\label{98122}
\end{align}

Mass($L^2$ norm)
\begin{align}
M(u)&:=\left(\int_{\mathbb{R}^N}
|u(x,t)|^2dx\right)^{\frac{1}{2}};\label{6193}
\end{align}

Energy
\begin{align}
E(u):=\frac{1}{2}\int_{\mathbb{R}^N}\left(|\nabla
u|^2+V(x)|u|^2\right)dx-\frac{1}{2}\int_{\mathbb{R}^N}F(x,|u|^2)dx-\frac{1}{4}\int_{\mathbb{R}^N}
(W\star|u|^2)|u|^2dx.\label{6194}
\end{align}

In \cite{Cazenave2}, Cazenave obtained some sufficient conditions on blowup and
global existence of the solution to (\ref{system0})
with (V1), (f1) and (W1). The following
two theorems can be looked as the parallel results to Corollary 6.1.2 and Theorem 6.5.4 of \cite{Cazenave2} respectively.

 {\bf Theorem 1.} (Global Existence) {\it Assume that $u_0 \in
\Sigma$, $(V2)$ and $(f1)$ are true, and
\begin{align}
W^+\in
L^q(\mathbb{R}^N)+L^{\infty}(\mathbb{R}^N)\label{98124}
\end{align}
for some $q\geq 1$, $q\geq \frac{N}{2}$(and $q>1$
if $N=2$). Here $W^+=\max(W,0)$. Suppose further that there exist
constants $c_1$ and $c_2$ such that $F(x,|u|^2)\leq c_1|u|^2+c_2|u|^{2p+2}$
with $0<p<\frac{2}{N}$. Then the solution of
(\ref{system0}) exists globally. That is,
$$
\|u(\cdot,t)\|_{\Sigma}<+\infty \quad {\rm for \ all} \
0<t<+\infty.
$$
}

{\bf Theorem 2.} (Blowup in Finite Time) {\it Assume that $u_0\in \Sigma$, $|x|u_0\in
L^2(\mathbb{R}^N)$, $(V2)$,
$(f1)$ and $(W1)$ are true. Suppose further that
 \begin{align}
& (N+2)F(x,|u|^2)-N|u|^2f(x,|u|^2)\leq 0,\label{9821x1}\\
& 2V(x)+(x\cdot \nabla V)\geq 0\quad
a.e.,\label{9821x2}\\
& 2W(x)+(x\cdot \nabla W)\leq 0\quad a.e.\label{9821x3}
\end{align}
If \quad (1) $E(u_0)<0$ or
(2) $E(u_0)=0$  and $\Im\int_{\mathbb{R}^N}(x\cdot\nabla
u_0)\bar{u}_0dx<0$,\\
then the solution of (\ref{system0}) will blow up in finite time.
That is, there exists $T_{\max}<\infty$ such that $$
\lim_{t\rightarrow T_{\max}}\|u(\cdot,t)\|_{\Sigma}=\infty.$$}

Denote
\begin{align}
Q(u)&:=2\int_{\mathbb{R}^N}|\nabla
u|^2dx-\int_{\mathbb{R}^N}(x\cdot \nabla
V)|u|^2dx\nonumber\\
&\quad+N\int_{\mathbb{R}^N}[F(x,|u|^2)-|u|^2f(x,|u|^2)]dx+\frac{1}{2}\int_{\mathbb{R}^N}((x\cdot
\nabla W)\star |u|^2)|u|^2dx.\label{963xw1}
\end{align}

We will establish the first type of sharp threshold as follows.

{\bf Theorem 3.} (Sharp Threshold I) {\it Assume that $V(x)\equiv
0$ and $W\in L^q(\mathbb{R}^N)$ with $\frac{N}{4}<q<\frac{N}{2}$.
Suppose further that $f(x,0)=0$ and there exist constants $c_1, c_2,
c_3>0$ and $\frac{2}{N}<p_1, p_2,l<\frac{2}{(N-2)^+}$ such that
\begin{align}
lF(x,|u|^2)\leq |u|^2f(x,|u|^2)-F(x,|u|^2)\leq c_1|u|^{2p_1+2}+c_2|u|^{2p_2+2},\label{97251}\\
NlW(x)+(x\cdot \nabla W)\leq 0\leq c_3W(x)+(x\cdot \nabla
W).\label{9828xw2}
\end{align}
Let $\omega$ be a positive constant satisfying
\begin{align} d_I:=\inf_{\{u\in \Sigma\setminus \{0\};
Q(u)=0\}}\left(\omega\|u\|_2^2+E(u)\right)>0,\label{9651}\end{align}
where $Q(u)$ is defined by (\ref{963xw1}). Suppose that $u_0\in
H^1(\mathbb{R}^N)$ satisfies
$$\omega\|u_0\|_2^2+E(u_0)<d_I.$$
Then

(1). If $Q(u_0)>0$, the solution of (\ref{system0}) exists
globally;

(2). If $Q(u_0)<0$, $|x|u_0\in L^2(\mathbb{R}^N)$ and $\Im
\int_{\mathbb{R}^N}(x\cdot \nabla u_0)\bar{u}_0dx<0$, the solution
of (\ref{system0}) blows up in finite time.}

{\bf Remark 1.1.} Theorem 3 is only suitable for (\ref{system0})
with $V(x)\equiv 0$. To establish the sharp threshold for
(\ref{system0}) with $V(x)\neq 0$, we will construct a type of
cross constrained variational problem and establish some
cross-invariant manifolds. First, we introduce some functionals as
follows:
\begin{align}
&I_{\omega}(u)=\omega
\|u\|_2^2+E(u),\label{8892}\\
&S_{\omega}(u)=2\omega \|u\|^2_2+\int_{\mathbb{R}^N}\left\{|\nabla
u|^2+V(x)|u|^2-f(x,|u|^2)|u|^2-(W\star|u|^2)|u|^2\right\}dx.\label{8181}
\end{align}
Denote the Nehari manifold
\begin{align}
\mathcal{N}&:=\{u\in \Sigma\setminus\{0\}, \
S_{\omega}(u)=0\},\label{9630x1}\end{align} and cross-manifold
\begin{align}
 \mathcal{CM}&:=\{u\in \Sigma\setminus\{0\},
\ S_{\omega}(u)<0,\ Q(u)=0\}.\label{5307}
\end{align}
And define
\begin{align}
d_{\mathcal{N}}&:=\inf_{\mathcal{N}} I_{\omega}(u),\label{5283}\\
d_{\mathcal{M}}&:=\inf_{\mathcal{CM}}I_{\omega}(u),\label{5308}\\
d_{II}&:= \min (d_{\mathcal{N}},\
d_{\mathcal{M}}).\label{9616w1}\end{align}

In Section 5, we will prove that $d_{II}$ is always positive.
Therefore, it is reasonable to define the following cross-manifold
\begin{align}
\mathcal{K}:&=\{u\in \Sigma\setminus\{0\}: I_{\omega}(u)<d_{II}, \
S_{\omega}(u)<0, \ Q(u)<0\}. \label{5316}
\end{align}

We give the second type of sharp threshold as follows

{\bf Theorem 4.} (Sharp Threshold II) {\it Assume that (f1), (W1)
and (\ref{97251}). Suppose that
\begin{align}
W(x)\geq 0, \quad NlW(x)+(x\cdot \nabla W)\leq 0\label{9816x1}
\end{align} and there exists a positive
constant $c$ such that
\begin{align}
NlV(x)+(x\cdot \nabla V)\geq cV(x)\geq 0\label{97252}
\end{align}
with the same $l$ in (\ref{97251}). Assume further that the
function $f(x,|u|^2)$ satisfies $f(x,0)=0$ and
\begin{align}
&f(x,|u|^2)\leq f(x,k^2|u|^2),\quad
f'_s(x,k^2|u|^2)\leq f'_s(x,|u|^2),\label{97253}\\
&F(x,k^2|u|^2)-k^2|u|^2f(x,k^2|u|^2)\leq k^2[F(x,|u|^2)-|u|^2f(x,|u|^2)]\label{9828w1}
\end{align}
for $k>1$. Here $f'_s(x,z)$ is the value of the partial derivative of $f(x,s)$ with respect to $s$ at the point $(x,z)$.  If $u_0\in \Sigma$ satisfies $|x|u_0\in
L^2(\mathbb{R}^N)$ and  $I_{\omega}(u_0)=\omega
\|u_0\|_2^2+E(u_0)<d_{II}$, then the solution of (\ref{system0})
blows up in finite time if and only if $u_0\in \mathcal{K}$.}

{\bf Remark 1.2.} (1) $f(x,|u|^2)\leq f(x,k^2|u|^2)$ implies that
$k^2F(x,|u|^2)\leq F(x,k^2|u|^2)$ for $k>1$.

(2) The blowup of solution to (\ref{system0}) will benefit from the condition $V(x)\geq 0$.
In some cases, the blowup of the solution to (\ref{system0})  can be
delayed or prevented by the introduction of potential(see
\cite{Carles} and the references therein).

This paper is organized as follows: In Section 2, we will prove Proposition 1.1, recall
some results of \cite{Cazenave2} and
give some other properties. In Section 3, we will prove
Theorem 1 and 2. In Section 4, we establish
the sharp threshold for (\ref{system0}) with $V(x)\equiv 0$. In Section 5, we will prove Theorem 4.

\section{Preliminaries}
\qquad In the sequel, we use $C$ and $c$ to denote various finite
constants, their exact values may vary from line to line.

First, we will give the proof of Proposition 1.1.

{\bf The proof of Proposition 1.1:}
If (V1) is true, then there exist $V_1(x)\in
L^r(\mathbb{R}^N)$ with $r\geq 1$,
$r>\frac{N}{2}$, and $V_2(x)\in L^{\infty}(\mathbb{R}^N)$ such that
$$V(x)=V_1(x)+V_2(x).$$
Noticing that $0<\frac{2r}{r-1}<\frac{2N}{N-2}$, using H\"{o}lder's and Soblev's inequalities, we have
\begin{align}
&\quad\int_{\mathbb{R}^N} V(x)|u|^2dx=\int_{\mathbb{R}^N} V_1(x)|u|^2dx+\int_{\mathbb{R}^N} V_2(x)|u|^2dx\nonumber\\
&\leq \left(\int_{\mathbb{R}^N} |V(x)|^rdx\right)^{\frac{1}{r}}
\left(\int_{\mathbb{R}^N} |u|^{\frac{2r}{r-1}}dx\right)^{\frac{r-1}{r}}+C\int_{\mathbb{R}^N} |u|^2dx\nonumber\\
&\leq C\int_{\mathbb{R}^N} |\nabla u|^2dx+C\int_{\mathbb{R}^N} |u|^2dx\label{104111}
\end{align}
for any $u\in H^1(\mathbb{R}^N)$. Consequently, we have
$$
\|u\|_{H^1}\leq \|u\|_{\Sigma} \leq C\|u\|_{H^1},
$$
which means that $\Sigma=H^1(\mathbb{R}^N)$ if
$V\in L^r(\mathbb{R}^N)+L^{\infty}(\mathbb{R}^N)$ for $r\geq
1$, $r>\frac{N}{2}$. By the results of  Theorem 3.3.1 in \cite{Cazenave2},
we have the local well-posedness result of (\ref{system0}) in $H^1(\mathbb{R}^N)$.

If (V2), (f1) and (W1) are true, similar to the proof of Theorem 3.5 in \cite{Oh}, we can
establish the local well-posedness result of (\ref{system0}) in $\Sigma$. Roughly, we only need to replace
$|u|^{p+1}u$ by $f(x,|u|^2)u+(W\star|u|^2)u$ in the proof, and we can obtain the similar results under the assumptions of (V2), (f1) and (W1). We omit the detail here.\hfill $\Box$

Noticing that $\Im h(u)\bar{u}=0$ and $h(u)=H'(u)$, following the
method of \cite{Glassey} and the discussion in Chapter 3 of
\cite{Cazenave2}, one can obtain the conservation of mass
and energy. We give the following proposition without proof.

{\bf Proposition 2.1.} {\it Assume that $u(x,t)$ is a solution of
(\ref{system0}). Then
\begin{align*}
M(u)&=\left(\int_{\mathbb{R}^N}
|u(x,t)|^2dx\right)^{\frac{1}{2}}=\left(\int_{\mathbb{R}^N}
|u_0(x)|^2dx\right)^{\frac{1}{2}}=M(u_0),
\\
E(u)&=\frac{1}{2}\int_{\mathbb{R}^N}\left\{|\nabla
u|^2+V(x)|u|^2-F(x,|u|^2)\right\}dx-G(|u|^2)=E(u_0)
\end{align*}}
for any $0\leq t<T_{\max}$.

We will recall some results on blowup and global
existence of the solution to (\ref{system0}) with (V1), (f1) and (W1).

{\bf Theorem A} (Corollary 6.1.2 of \cite{Cazenave2}) {\it Assume
that $(V1)$, $(f1)$ and (\ref{98124}). Suppose that there
exist $A\geq 0$ and $0\leq p<\frac{2}{N}$ such that
\begin{align}
F(|u|^2)\leq A|u|^2(1+|u|^{2p}).\label{98123}
\end{align}
 Then the maximal strong $H^1$-solution of
(\ref{system0}) is global and $\sup\{\|u\|_{H^1}:t\in
\mathbb{R}\}<\infty$ for every $u_0\in H^1(\mathbb{R}^N)$.
 }

{\bf Theorem B}  (Theorem 6.5.4 of \cite{Cazenave2})  {\it  Assume
that $(V1)$, $(f1)$, $(W1)$ and
(\ref{9821x1})--(\ref{9821x3}). If $u_0\in
H^1(\mathbb{R}^N)$, $|x|u_0\in L^2(\mathbb{R}^N)$ and
$E(u_0)<0$, then the $H^1$-solution of (\ref{system0}) will blow
up in finite time.
 }

Let $J(t)=\int_{\mathbb{R}^N} |x|^2|u|^2dx$. After some elementary
computations, we obtain
$$ J'(t)=4\Im \int_{\mathbb{R}^N}\{(x\cdot \nabla u) \bar{u}dx,
\quad J''(t)=4Q(u).
$$

We have the following proposition

 {\bf Proposition 2.2.}  {\it
Assume that $u(x,t)$ is a solution of (\ref{system0}) with $u_0\in \Sigma$ and $|x|u_0\in
L^2(\mathbb{R}^N)$. Then the solution to
(\ref{system0}) will blow up in finite time if either

(1) there exists a constant $c<0$ such that $J''(t)=4Q(u)\leq c<0$
or

(2) $J''(t)=4Q(u)\leq 0$ and $J'(0)=\Im \int_{\mathbb{R}^N}(x\cdot
\nabla u_0)\bar{u}_0dx<0$.}

{\bf Proof:} Since $u_0\in \Sigma$ and $|x|u_0\in
L^2(\mathbb{R}^N)$, we have
$$
|J'(0)|<4\int_{\mathbb{R}^N}|(x\bar{u}_0||
\nabla u_0)|dx\leq 8\int_{\mathbb{R}^N} (|\nabla u_{10}|^2+|xu_{10}|^2)dx<+\infty.
$$

(1) If $J''(t)\leq c<0$, integrating it from $0$ to
$t$, we get $J'(t)<ct+J'(0)$. Since $c<0$,  we know that
there exists a $t_0\geq \max(0,\frac{J'(0)}{-c})$ such that
$J'(t)<J'(t_0)<0$ for $t>t_0$. On the other hand, we have
\begin{align}
0\leq
J(t)=J(t_0)+\int_{t_0}^tJ'(s)ds<J(t_0)+J'(t_0)(t-t_0),\label{95312}
\end{align}
which implies that there exists a $T_{\max}<+\infty$ satisfying
\begin{align}
\lim_{t\rightarrow
T_{\max}}J(t)=0.\label{2010781}
\end{align}
Using the inequality
\begin{align}
\|g\|^2_2\leq \frac{2}{N}\|\nabla g\|_2\|xg\|_2 \quad {\rm if }\ g\in H^1(\mathbb{R}^N, \ xg\in L^2(\mathbb{R}^N) \label{2010782}
\end{align}
and noticing that $\|u(\cdot,t)\|_2=\|u_0\|_2$, we have
$$
\lim_{t\rightarrow
T_{\max}}\|u\|_{\Sigma}=+\infty.
$$

(2) Similar to (\ref{95312}), we can get
$$0\leq
J(t)\leq J(0)+J'(0)t,$$ which implies that the solution will blow
up in a finite time $T_{\max}\leq\frac{J(0)}{-J'(0)}$. \hfill
$\Box$

\section{The sufficient conditions on global existence and blowup in finite time}
\qquad In this section, we will prove Theorem 1 and 2, which give
some sufficient conditions on global existence and blowup of
the solution to (\ref{system0}).

{\bf The proof of Theorem 1:} Letting $W^+=W_1+W_2$, where $W_1\in
L^{\infty}$ and $W_2\in L^q$ with $q>\frac{N}{2}$,
using H\"{o}lder's and Young's inequalities, we obtain
$$
\int_{\mathbb{R}^N}(W_2\star(uv))wzdx\leq \|W_2\|_{L^q}\|u\|_{L^r}\|v\|_{L^r}\|w\|_{L^r}\|z\|_{L^r}
$$
with $r=\frac{4q}{2q-1}$. Especially, we have
\begin{align}
\int_{\mathbb{R}^N}(W_2\star|u|^2)|u|^2dx\leq \|W_2\|_{L^q}\|u\|^4_{L^r}.\label{20107101}
\end{align}
Using (\ref{20107101}) and
Gagliardo-Nirenberg's inequality, we get
\begin{align}
\frac{1}{4}\int_{\mathbb{R}^N} (W\star|u|^2)|u|^2dx&\leq
\|W_1\|_{L^{\infty}}\|u\|_{L^2}^4+\|W_2\|_{L^q}\|u\|_{L^{\frac{4q}{2q-1}}}^4\nonumber\\
&\leq\|W_1\|_{L^{\infty}}\|u\|_{L^2}^4+C\|W_2\|_{L^q}\| \nabla
u\|_{L^2}^{\frac{N}{q}}\|u\|_{L^2}^{\frac{4q-N}{q}}.\label{98141}
\end{align}
Using Young's inequality, from (\ref{98141}), we have
\begin{align}
C\|W_2\|_{L^q}\|\nabla
u\|_{L^2}^{\frac{N}{q}}\|u\|_{L^2}^{\frac{4q-N}{q}}\leq
\varepsilon \|\nabla u\|_{L^2}^2+C(\varepsilon,
\|W_2\|_{L^q})\|u\|_{L^2}^{\frac{8q-2N}{2q-N}}\label{98142}
\end{align}
for some $\varepsilon>0$. Noticing that $F(x,|u|^2)\leq c_1|u|^2+c_2|u|^{2p+2}$,
using Gagliardo-Nirenberg's inequality
and (\ref{98142}) with $\varepsilon=\frac{1}{4}$, we get
\begin{align}
E(u_0)&=\frac{1}{2}\left(\int_{\mathbb{R}^N}\left\{|\nabla
u_0|^2+V(x)|u_0|^2-F(x,|u_0|^2)\right\}dx\right)-\frac{1}{4}\int_{\mathbb{R}^N}
(W\star|u_0|^2)|u_0|^2dx\nonumber\\
&=\frac{1}{2}\left(\int_{\mathbb{R}^N}\left\{|\nabla
u|^2+V(x)|u|^2-F(x,|u|^2)\right\}dx\right)-\frac{1}{4}\int_{\mathbb{R}^N}
(W\star|u|^2)|u|^2dx\nonumber\\
&\geq \frac{1}{2}\left(\int_{\mathbb{R}^N}\left\{|\nabla
u|^2+V(x)|u|^2-c_1|u|^2-c_2|u|^{2p+2}\right\}dx\right)\nonumber\\
&\qquad-\|W_1\|_{L^{\infty}}\|u\|_{L^2}^4-C\|W_2\|_{L^q}\|\nabla
u\|_{L^2}^{\frac{N}{q}}\|u\|_{L^2}^{\frac{4q-N}{q}}\nonumber\\
&\geq \frac{1}{2}\left(\int_{\mathbb{R}^N}\left\{|\nabla
u|^2+V(x)|u|^2-c_1|u|^2\right\}dx\right)\nonumber\\
&\qquad-c_2C_N\left(\int_{\mathbb{R}^N}|\nabla
u|^2dx\right)^{\frac{pN}{2}}\left(\int_{\mathbb{R}^N}|
u|^2dx\right)^{\frac{2+p(2-N)}{2}}\nonumber\\
&\qquad -\|W_1\|_{L^{\infty}}\|u\|_{L^2}^4-\frac{1}{4} \|\nabla
u\|_{L^2}^2-C\|u\|_{L^2}^{\frac{8-2N}{2q-N}}.\label{97161}
\end{align}
Since $\|u\|_2=\|u_0\|_2$, from (\ref{97161}), we can obtain
\begin{align}
&\quad 4E(u_0)+C\|u_0\|_{L^2}^2+C\|u_0\|_{L^2}^4+C\|u_0\|_{L^2}^{\frac{8-2N}{2q-N}}\nonumber\\
&\geq \int_{\mathbb{R}^N}V(x)|u|^2dx+\int_{\mathbb{R}^N}|\nabla
u|^2dx\left(1-c\left\{\int_{\mathbb{R}^N}|\nabla
u|^2dx\right\}^{\frac{pN}{2}-1}\right).\label{97162}
\end{align}
Since $p<\frac{2}{N}$ means that $\frac{pN}{2}-1<0$, (\ref{97162})
implies that $\|u\|^2_{\Sigma}$ is always controlled by
$4E(u_0)+C\|u_0\|_{L^2}^2+C\|u_0\|_{L^2}^4+C\|u_0\|_{L^2}^{\frac{8-2N}{2q-N}}$.
That is, the solution of (\ref{system0}) exists globally.\hfill
$\Box$

{\bf Remark 3.1.} We will give some examples of $V(x)$, $f(x,|u|^2)$ and
$W(x)$. It is easy to verify that they satisfy the conditions of Theorem 1.

Example 1. $V(x)=|x|^2$, $W(x)=e^{-\pi|x|^2}$ and $f(x,|u|^2)=b|u|^{2p}$ with $b$ is a real constant and $0<p<\frac{2}{N}$.

Example 2. $V(x)=|x|^2$, $W(x)=\frac{|x|^2}{1+|x|^2}$ and
$f(x,|u|^2)=b|u|^{2p}\ln(1+|u|^2)$ with $b$ is a real constant and
$0<p<\frac{2}{N}$.

{\bf The proof of Theorem 2:} Set \begin{align} y(t)=J'(t)=4\Im
\int_{\mathbb{R}^N}(x\cdot \nabla u)
\bar{u}dx.\label{9741}\end{align}
Using (\ref{9821x1})-(\ref{9821x3}), we have
\begin{align}
y'(t)&=8\int_{\mathbb{R}^N}|\nabla
u|^2dx-4\int_{\mathbb{R}^N}(x\cdot \nabla
V)|u|^2dx\nonumber\\
&\quad
+4N\int_{\mathbb{R}^N}[F(x,|u|^2)-|u|^2f(x,|u|^2)]dx
+2\int_{\mathbb{R}^N}\left\{(x\cdot
\nabla W)\star
|u|^2\right\}|u|^2dx\nonumber\\
&=16E(u)+4\int_{\mathbb{R}^N}\left([-2V(x)-(x\cdot \nabla
V)]|u|^2+[(N+2)F(x,|u|^2)-N|u|^2f(x,|u|^2)]\right)dx\nonumber\\
&\quad+2\int_{\mathbb{R}^N}[\{2W+(x\cdot \nabla W)\}\star
|u|^2]|u|^2dx\leq 16E(u)=16E(u_0)<0.\label{9742}\end{align}
From (\ref{9741}) and (\ref{9742}), we obtain
\begin{align}
\|xu(x,t)\|_{L^2}^2\leq
\|xu_0\|_{L^2}^2+4t\Im\int_{\mathbb{R}^N}\bar{u}_0(x\cdot\nabla
u_0)dx+8t^2E(u_0).\label{9814x1}
\end{align}
Since $\|xu(x,t)\|_{L^2}^2\geq 0$, whether (1)
or (2), (\ref{9814x1}) will be absurd for $t>0$ large
enough. Therefore, the solution of (\ref{system0}) will blow up in
finite time.\hfill $\Box$

{\bf Remark 3.2.} We will give some examples of $V(x)$, $W(x)$ and
$f(x,|u|^2)$. It is easy to verify that they satisfy the conditions of Theorem 2.

Example 1. $V(x)=|x|^2$, $W(x)=|x|^{-2}$ and $f(x,|u|^2)=b|u|^{2p}$ with $b>0$ and $p>\frac{2}{N}$ with $N\geq 3$.

Example 2. $V(x)=|x|^2$, $W(x)=|x|^{-2}$ and $f(x,|u|^2)=b|u|^{2p}\ln(1+|u|^2)$
with $b>0$ and $p\geq \frac{2}{N}$ with $N\geq 3$.

\section{The sharp threshold for global existence and blowup of the solution
to (\ref{system0}) with $V(x)\equiv 0$ and $ W\in
L^q(\mathbb{R}^N)$  with $\frac{N}{4}<q<\frac{N}{2}$}

\qquad In this section, we will establish the sharp threshold for
global existence and blowup of the solution to (\ref{system0})
with $V(x)\equiv 0$ and $ W\in L^q(\mathbb{R}^N)$  with
$\frac{N}{4}<q<\frac{N}{2}$.

{\bf The proof of Theorem 3.} We will proceed in four steps.

Step 1. We will prove $d_I>0$. $u\in H^1(\mathbb{R}^N)\setminus
\{0\}$ and $Q(u)=0$ mean that
\begin{align*}
2\int_{\mathbb{R}^N}|\nabla
u|^2dx&=N\int_{\mathbb{R}^N}[|u|^2f(x,|u|^2)-F(x,|u|^2)]dx-\frac{1}{2}\int_{\mathbb{R}^N}
\{(x\cdot \nabla W)\star|u|^2\}|u|^2dx\nonumber\\
&\leq
\frac{N(l+1)}{l}\int_{\mathbb{R}^N}[c_1|u|^{2p_1+2}+c_2|u|^{2p_2+2}]dx+C\int_{\mathbb{R}^N}
(W\star|u|^2)|u|^2dx\nonumber\\
&\leq
C\|u\|_{2p_1+2}^{2p_1+2}+C\|u\|_{2p_2+2}^{2p_2+2}+C\|W\|_{L^q}\|u\|_{L^{\frac{4q}{2q-1}}}^4.
\end{align*}
Using Gagliardo-Nirenberg's and H\"{o}lder's inequalities, we can
get
\begin{align*}
2&\leq C(\|\nabla
u\|^2_2)^{\frac{Np_1}{2}}(\|u\|_2^2)^{p_1+1-\frac{Np_1}{2}}+C(\|\nabla
u\|^2_2)^{\frac{Np_2}{2}}(\|u\|_2^2)^{p_2+1-\frac{Np_2}{2}}\nonumber\\
&\quad+C(\|\nabla
u\|^2_2)^{\frac{N}{2q}}(\|u\|_2^2)^{\frac{4q-N}{2q}}\\
&\leq C\left\{\left(\|\nabla
u\|_2^2+\|u\|_2^2\right)^{p_1+1}+\left(\|\nabla
u\|_2^2+\|u\|_2^2\right)^{p_2+1}+(\|\nabla
u\|_2^2+\|u\|_2^2)^2\right\}.
\end{align*}
That is,
\begin{align}
\|\nabla u\|_2^2+\|u\|_2^2\geq C>0\label{96155}
\end{align}
if $Q(u)=0$ and $u\in H^1(\mathbb{R}^N)\setminus \{0\}$.

On the other hand, if $Q(u)=0$, we have
\begin{align*}
2\int_{\mathbb{R}^N}|\nabla
u|^2dx&=N\int_{\mathbb{R}^N}[|u|^2f(x,|u|^2)-F(x,|u|^2)]dx-\frac{1}{2}\int_{\mathbb{R}^N}
\{(x\cdot \nabla W)\star|u|^2\}|u|^2dx\nonumber\\
&\geq
Nl\int_{\mathbb{R}^N}F(x,|u|^2)dx+\frac{Nl}{2}\int_{\mathbb{R}^N}
\{W\star|u|^2\}|u|^2dx,
\end{align*}
that is,
\begin{align}
-\frac{1}{2}\int_{\mathbb{R}^N}F(x,|u|^2)dx-\frac{1}{4}\int_{\mathbb{R}^N}
\{W\star|u|^2\}|u|^2dx\geq -\frac{1}{Nl}\int_{\mathbb{R}^N}|\nabla
u|^2dx.\label{9821w1}
\end{align}
Using (\ref{9821w1}), we can obtain
\begin{align*}
\omega
\|u\|_2^2+E(u)&=\omega\|u\|_2^2+\frac{1}{2}\int_{\mathbb{R}^N}|\nabla
u|^2dx-\frac{1}{2}\int_{\mathbb{R}^N}F(x,|u|^2)dx-\frac{1}{4}\int_{\mathbb{R}^N}
\{W\star|u|^2\}|u|^2dx\\
&\geq
\omega\|u\|_2^2+(\frac{1}{2}-\frac{1}{Nl})\int_{\mathbb{R}^N}|\nabla
u|^2dx\\ &\geq \min\{\omega,
(\frac{1}{2}-\frac{1}{Nl})\}\left(\|\nabla
u\|_2^2+\|u\|_2^2\right)\geq C>0
\end{align*}
from (\ref{96155}). Hence $$d_I>0.$$

Step 2. Denote
$$
K_+=\{u\in H^1(\mathbb{R}^N)\setminus\{0\},\ Q(u)>0,\
\omega\|u\|_2^2+E(u)<d_I\}
$$
and
$$
K_-=\{u\in H^1(\mathbb{R}^N)\setminus\{0\},\ Q(u)<0,\
\omega\|u\|_2^2+E(u)<d_I\}.
$$
We will prove that $K_+$ and $K_-$ are invariant sets of
(\ref{system0}) with $V(x)\equiv 0$ and $ W\in L^q(\mathbb{R}^N)$
with $\frac{N}{4}<q<\frac{N}{2}$. That is, we need to show
that $u(\cdot,t)\in \mathcal{K}$ for all $t\in (0,T_{\max})$ if
$u_0\in K_+$. Since $\|u\|_2$ and $E(u)$ are conservation
quantities for (\ref{system0}), we have
\begin{align}
u(\cdot,t)\in H^1(\mathbb{R}^N)\setminus\{0\}, \quad
\omega\|u(\cdot,t)\|_2^2+E(u(\cdot,t))<d_I\label{9626z2}\end{align}
for all $t\in (0,T_{\max})$ if $u_0\in K_+$. We need to prove that
$Q(u(\cdot,t))>0$. Otherwise, assume that there exists a $t_1\in
(0,T_{\max})$ satisfying $Q(u(\cdot,t_1))=0$ by the continuity.
Note that (\ref{9626z2}) implies
$$\omega\|u(\cdot,t_1)\|_2^2+E(u(\cdot,t_1))<d_I.$$
However, the inequality above and $Q(u(\cdot,t_1))=0$ are
contradictions to the definition of $d_I$. Therefore,
$Q(u(\cdot,t))>0$. Consequently, (\ref{9626z2}) and
$Q(u(\cdot,t))>0$  imply that $u(\cdot,t)\in K_+$. That is,
$K_+$ is a invariant set of (\ref{system0}) with $V(x)\equiv 0$
and $ W\in L^q(\mathbb{R}^N)$  with $\frac{N}{4}<q<\frac{N}{2}$.
Similarly, we can prove that $K_-$ is also a invariant set of
(\ref{system0}) with $V(x)\equiv 0$ and $ W\in L^q(\mathbb{R}^N)$
with $\frac{N}{4}<q<\frac{N}{2}$.

Step 3. Assume that $Q(u_0)>0$ and $\omega\|u_0\|_2^2+E(u_0)<d_I$.
By the results of Step 2, we have $Q(u(\cdot,t))>0$ and
$\omega\|u(\cdot,t)\|_2^2+E(u(\cdot,t))<d_I$. That is,
\begin{align*}
-2\|\nabla
u(\cdot,t)\|_2^2&<-N\int_{\mathbb{R}^N}[|u|^2f(x,|u|^2)-F(x,|u|^2)]dx+\frac{1}{2}\int_{\mathbb{R}^N}
\{(x\cdot \nabla W)\star|u|^2\}|u|^2dx\\
&<-Nl\int_{\mathbb{R}^N}F(x,|u|^2)dx-\frac{Nl}{2}\int_{\mathbb{R}^N}
\{W\star|u|^2\}|u|^2dx,\end{align*} and \begin{align*}
d_I&>\omega\|u(\cdot,t)\|_2^2+\frac{1}{2}\|\nabla
u(\cdot,t)\|_2^2-\frac{1}{2}\int_{\mathbb{R}^N}F(x,|u|^2)dx-\frac{1}{4}\int_{\mathbb{R}^N}
\{W\star|u|^2\}|u|^2dx.
\end{align*}
The two inequalities imply that
\begin{align*}
\omega\|u(\cdot,t)\|_2^2+(\frac{1}{2}-\frac{1}{Nl})\|\nabla
u(\cdot,t)\|_2^2<d_I.
\end{align*} which means that
$$
\|u(\cdot,t)\|_{H^1(\mathbb{R}^N)}<\infty,
$$
i.e., the solution  exists globally.

Step 4. Assume that $Q(u_0)<0$ and $\omega\|u_0\|_2^2+E(u_0)<d_I$.
By the results of Step 2, we obtain $Q(u(\cdot,t))<0$ and
$\omega\|u(\cdot,t)\|_2^2+E(u(\cdot,t))<d_I$. Hence we get
$$J''(t)=4Q(u)<0,\quad J'(0)=4\Im \int_{\mathbb{R}^N}(x\cdot
\nabla u_0)\bar{u}_0dx<0.$$ By the results of Proposition 2.2, the
solution will blow up in finite time.\hfill $\Box$

As a corollary of Theorem 3, we obtain the sharp threshold for
global existence and blowup of the solution of (\ref{972310}) as
follows.

{\bf Corollary 4.1.} {\it Assume that $f(x,0)=0$ and (\ref{97251}).
Let $\omega$ be a positive constant satisfying
\begin{align} d'_I:=\inf_{\{u\in \Sigma\setminus \{0\};
Q_1(u)=0\}}\left(\omega\|u\|_2^2+E(u)\right)>0.\label{9651}\end{align}
Here
\begin{align}
Q_1(u)&:=2\int_{\mathbb{R}^N}|\nabla u|^2dx
+N\int_{\mathbb{R}^N}[F(x,|u|^2)-|u|^2f(x,|u|^2)]dx.\label{9828xw1}
\end{align}

Suppose that $u_0\in H^1(\mathbb{R}^N)$ satisfies
$$\omega\|u_0\|_2^2+E(u_0)<d'_I.$$

Then

(1). If $Q_1(u_0)>0$, the solution of (\ref{972310}) exists
globally;

(2). If $Q_1(u_0)<0$, $|x|u_0\in L^2(\mathbb{R}^N)$ and $\Im
\int_{\mathbb{R}^N}(x\cdot \nabla u_0)\bar{u}_0dx<0$, the solution
of (\ref{972310}) blows up in finite time.}

{\bf Remark 4.1.}  In Theorem 1.5 of \cite{Tao2}, Tao et al.
proved that:

{\it Assume that $u(x,t)$ is  a solution
of (\ref{972310}) with $f(x,|u|^2)u=\mu |u|^{p_1}u+\nu |u|^{p_2}u$, where $\mu>0$, $\nu>0$, $ \frac{4}{N}\leq
p_1<p_2\leq\frac{4}{N-2}$ with $N\geq 3$, $\Im \int_{\mathbb{R}^N}(x\cdot \nabla
u_0)\bar{u}_0dx<0$, $|x|u_0\in L^2(\mathbb{R}^N)$ and $E(u_0)<0$. Then blowup occurs.}

Corollary 4.1 improve the result above. In fact,
if $f(x,|u|^2)u=\mu |u|^{p_1}u+\nu |u|^{p_2}u$, then
\begin{align*}
Q_1(u)&=4E(u)-\frac{(Np_1-4)\mu
}{(p_1+2)}\|u\|_{p_1+2}^{p_1+2}-\frac{(Np_2-4)\nu
}{(p_2+2)}\|u\|_{p_2+2}^{p_2+2}\leq E(u),
\end{align*}
hence $E(u_0)<0$ implies that $Q_1(u_0)<0$. That is, our blowup
condition is weaker than theirs. On the other hand, our conclusion is still
true if $0<E(u_0)<d'_I-\omega \|u_0\|_2^2$ with $Q_1(u_0)<0$, $\Im \int_{\mathbb{R}^N}(x\cdot \nabla
u_0)\bar{u}_0dx<0$ and $|x|u_0\in L^2(\mathbb{R}^N)$. In
other words, our result is stronger than theirs if
$\omega\|u_0\|_2^2+E(u_0)<d'_I$ with $Q_1(u_0)<0$, $\Im \int_{\mathbb{R}^N}(x\cdot \nabla
u_0)\bar{u}_0dx<0$ and $|x|u_0\in L^2(\mathbb{R}^N)$.

{\bf Remark 4.2.} We will give some examples of  $f(x,|u|^2)$ and $W(x)$.
It is easy to verify that they satisfy the conditions of Theorem 3.

Example 4.1.  $W(x)\equiv 0$, $f(x,|u|^2)=c|u|^{2q_1}+d|u|^{2q_2}$
with $c<0$, $d>0$ and $q_2>\frac{2}{N}$, $q_2>q_1>0$.

Example 4.2.  $W(x)\equiv 0$, $f(x,|u|^2)=b|u|^{2p}\ln(1+|u|^2)$
with $b>0$ and $p>\frac{2}{N}$.

Example 4.3.  Let $f(x,|u|^2)$ be one of those in Examples 4.1 and
4.2. And Let
\begin{align*} W(x)=\left\{\begin{array}{ll}&\frac{1}{|x|^{Nl}},\quad |x|\leq 1,\\
&\varphi(x),\quad  1\leq |x|\leq 2,\\
&\frac{1}{|x|^{K}}, \quad |x|\geq 2,
\end{array}
\right. \end{align*} where $2<Nl<\frac{N}{q}<K$, and $\varphi(x)$
satisfies
$$ Nl\varphi(x)+(x\cdot \nabla \varphi)\leq 0\leq
c_3\varphi(x)+(x\cdot \nabla \varphi)$$ when $1\leq |x|\leq 2$ and
makes $W(x)$ be smooth. Obviously, $W\in L^q(\mathbb{R}^N)$.

\section{Sharp threshold for the blowup and global existence of the solution to (\ref{system0})}

\qquad Theorem 4 extend the results of \cite{Zhang} to more general case. Moreover, we need subtle estimates and
more sophisticated analysis in the proof.

\subsection{Some invariant manifolds}
\qquad In this subsection, we will prove that $d_{\mathcal{N}},
d_{\mathcal{M}}, d_{II}>0$, and construct some invariant
manifolds.

{\bf Proposition 5.1.1.} {\it Assume that the conditions of Theorem
4 hold. Then $d_{\mathcal{N}}>0$.}

{\bf Proof:}  Assume that $u\in \Sigma\setminus \{0\}$ satisfying
$S_{\omega}(u)= 0$. Using Gagliardo-Nirenberg's and
Young's inequalities, we have
\begin{align}
&\qquad 2\omega\|u\|_2^2+\int_{\mathbb{R}^N}[|\nabla
u|^2+V(x)|u|^2]dx\nonumber\\
&=\int_{\mathbb{R}^N}|u|^2f(x,|u|^2)dx+\int_{\mathbb{R}^N}(W\star|u|^2)|u|^2dx\nonumber\\
&\leq\frac{l+1}{l}\int_{\mathbb{R}^N}[c_1|u|^{2p_1+2}+c_2|u|^{2p_2+2}]dx+\|W_1\|_{L^{\infty}}\|u\|_2^4
+\|W_2\|_{L^q}\|u\|_{L^{\frac{4q}{2q-1}}}^4\nonumber\\
&\leq C(\|\nabla
u\|^2_2)^{\frac{Np_1}{2}}(\|u\|_2^2)^{p_1+1-\frac{Np_1}{2}}+C(\|\nabla
u\|^2_2)^{\frac{Np_2}{2}}(\|u\|_2^2)^{p_2+1-\frac{Np_2}{2}}\nonumber\\
&\qquad+\|W_1\|_{L^{\infty}}\|u\|_2^4+C\|W_2\|_{L^q}\|
\nabla u\|_2^{\frac{N}{q}}\|u\|_2^{\frac{4q-N}{q}}\nonumber\\
&\leq C(\|\nabla
u\|^2_2)^{\frac{Np_1}{2}}(\|u\|_2^2)^{p_1+1-\frac{Np_1}{2}}+C(\|\nabla
u\|^2_2)^{\frac{Np_2}{2}}(\|u\|_2^2)^{p_2+1-\frac{Np_2}{2}}\nonumber\\
&\qquad+C\|u\|_2^4+\|\nabla
u\|_2^4+C(\|W_2\|_{L^q})\|u\|_2^4.\label{9814w1}\end{align}
Using H\"{o}lder's inequality,  from (\ref{9814w1}), we can obtain
\begin{align}
 &\qquad 2\omega\|u\|_2^2+\int_{\mathbb{R}^N}[|\nabla
u|^2+V(x)|u|^2]dx\nonumber\\
&\leq
C\left(2\omega\|u\|_2^2+\int_{\mathbb{R}^N}[|\nabla
u|^2+V(x)|u|^2]dx\right)^{p_1+1}\nonumber\\
&\quad +C\left(2\omega\|u\|_2^2+\int_{\mathbb{R}^N}[|\nabla
u|^2+V(x)|u|^2]dx\right)^{p_2+1}\nonumber\\
&\quad + C\left(2\omega\|u\|_2^2+\int_{\mathbb{R}^N}[|\nabla
u|^2+V(x)|u|^2]dx\right)^2.\label{97171}
\end{align}
(\ref{97171}) implies that
\begin{align}
2\omega\|u\|_2^2+\int_{\mathbb{R}^N}[|\nabla u|^2+V(x)|u|^2]dx\geq
C>0\label{97172}
\end{align}
for some positive constant $C$.

On the other hand, if $S_{\omega}(u)=0$, we get
\begin{align}
&\quad\omega \|u\|_2^2+\frac{1}{2}\int_{\mathbb{R}^N}(|\nabla
u|^2+V(x)|u|^2)dx\nonumber\\
&=\frac{1}{2}\int_{\mathbb{R}^N}f(x,|u|^2)|u|^2dx+\frac{1}{2}\int_{\mathbb{R}^N}(W\star|u|^2)|u|^2dx\nonumber\\
&\geq \min(l+1,2)\left(\frac{1}{2}\int_{\mathbb{R}^N}F(x,|u|^2)dx
+\frac{1}{4}\int_{\mathbb{R}^N}(W\star|u|^2)|u|^2dx\right).\label{9822x1}
\end{align}
From (\ref{9822x1}), we obtain
\begin{align}
I_{\omega}(u)&=\omega\|u\|_2^2+\frac{1}{2}\int_{\mathbb{R}^N}[|\nabla
u|^2+V(x)|u|^2-F(x,|u|^2)]dx-G(|u|^2)\nonumber\\
&\geq
\min\left(\frac{l}{2(l+1)},\frac{1}{4}\right)\left(2\omega\|u\|_2^2+\int_{\mathbb{R}^N}[|\nabla
u|^2+V(x)|u|^2]dx\right)\nonumber\\
&\geq C>0.\label{97173}
\end{align}

Consequently, $$ \qquad \qquad \qquad \qquad \qquad\qquad
d_{\mathcal{N}}=\inf_{\mathcal{N}}I_{\omega}(u)>C>0. \qquad \qquad
\qquad \qquad \qquad  \qquad \qquad
 \Box$$

Now, we will give some properties of $I_{\omega}(u)$,
$S_{\omega}(u)$ and $Q(u)$. We have a proposition as follows.

{\bf Proposition 5.1.2.} {\it Assume that
 $Q(u)$ and $S_{\omega}(u)$ are defined by (\ref{963xw1}) and
 (\ref{8181}). Then we have

(i) There at least exists a $w^{\star}\in \Sigma\setminus \{0\}$
such that
\begin{align}
S_{\omega}(w^{\star})=0,\quad Q(w^{\star})=0.\label{7181}
\end{align}

(ii) There at least exists a $u^*\in \Sigma\setminus \{0\}$  such
that
\begin{align}
 S_{\omega}(u^*)<0, \quad
Q(u^*)=0.\label{5306}
\end{align}}

{\bf Proof:} (i)  Noticing the assumptions on $V(x)$, $W(x)$ and
$f(x,|u|^2)$, similar to the proof of Theorem 1.7 in \cite{Rabinowitz},
it is easy to prove that there exists a $w^{\star}\in
\Sigma\setminus \{0\}$ satisfying
\begin{align}
2\omega w^{\star}+V(x)w^{\star}-\Delta
w^{\star}=f(x,|w^{\star}|^2)w^{\star}+(W\star|w^{\star}|^2)w^{\star}
\quad {\rm in} \ \mathbb{R}^N.\label{88101}
\end{align}
Multiplying (\ref{88101}) by $w^{\star}$ and integrating over
$\mathbb{R}^N$ by part, we can get $S_{\omega}(w^{\star})=0$.

Multiplying (\ref{88101}) by $(x\cdot \nabla w^{\star})$ and
integrating over $\mathbb{R}^N$ by part, we obtain the Pohozaev's
identity:
\begin{align}
&\quad
N\omega\|w^{\star}\|^2_2+\frac{N-2}{2}\int_{\mathbb{R}^N}|\nabla
w^{\star}|^2dx+\frac{N}{2}\int_{\mathbb{R}^N}V(x)|w^{\star}|^2dx+\frac{1}{2}\int_{\mathbb{R}^N}(x\cdot
\nabla V)|w^{\star}|^2dx\nonumber\\
&=\frac{N}{2}\int_{\mathbb{R}^N}F(x,|w^{\star}|^2)dx+\frac{N}{2}\int_{\mathbb{R}^N}(W\star|w^{\star}|^2)|w^{\star}|^2dx
+\frac{1}{2}\int_{\mathbb{R}^N}\{(x\cdot \nabla
W)\star|w^{\star}|^2\}|w^{\star}|^2dx.\label{88102}
\end{align}
From $S_{\omega}(w^{\star})=0$ and (\ref{88102}), we can get
$Q(w^{\star})=0$.

(ii) Letting $v_{k,\lambda}(x)=kw^{\star}(\lambda x)$ for $k>0$
and $\lambda>0$, we can obtain
\begin{align}
S_{\omega}(v_{k,\lambda})&=2\omega
k^2\int_{\mathbb{R}^N}|w^{\star}(\lambda
x)|^2dx+k^2\int_{\mathbb{R}^N}|\nabla
w^{\star}(\lambda x)|^2dx+k^2\int_{\mathbb{R}^N}V(x)|w^{\star}(\lambda x)|^2dx\nonumber\\
&\quad
-k^2\int_{\mathbb{R}^N}|w^{\star}(\lambda x)|^2f(x,k^2|w^{\star}(\lambda x)|^2)dx
-k^4\int_{\mathbb{R}^N}\left(W\star|w^{\star}(\lambda x)|^2\right)|w^{\star}(\lambda x)|^2dx,\label{95221}\\
Q(v_{k,\lambda})&=2k^2\int_{\mathbb{R}^N}|\nabla w^{\star}(\lambda
x)|^2dx-k^2\int_{\mathbb{R}^N}(x\cdot \nabla V)|w^{\star}(\lambda
x)|^2dx\nonumber\\&\qquad-N\int_{\mathbb{R}^N}[k^2|w^{\star}(\lambda
x)|^2f(x,k^2|w^{\star}(\lambda x)|^2)-F(x,k^2|w^{\star}(\lambda
x)|^2)dx\nonumber\\
&\qquad+\frac{k^4}{2}\int_{\mathbb{R}^N}\left((x\cdot\nabla
W)\star |w^{\star}(\lambda x)|^2\right)|w^{\star}(\lambda
x)|^2dx.\label{95222}
\end{align}
Looking $S_{\omega}(v_{k,\lambda})$ and $Q(v_{k,\lambda})$ as the
functions of $(k,\lambda)$, setting
$g(k,\lambda)=S_{\omega}(v_{k,\lambda})$ and
$\eta(k,\lambda)=Q(v_{k,\lambda})$, we get that $g(1,1)=0$ and
$\eta(1,1)=0$. And we want to prove that there exists a pair of $(k,\lambda)$
such that $g(k,\lambda)=S_{\omega}(v_{k,\lambda})<0$  and
$\eta(k,\lambda)=Q(v_{k,\lambda})=0$. Since $\eta(1,1)=0$, we know
that the image of $\eta(k,\lambda)$ and the plane $\eta=0$
intersect in the space of $(k,\lambda, \eta)$ and form a curve
$\eta(k,\lambda)=0$. Hence there exist many positive real number
pairs $(k,\lambda)$ relying on $w^{\star}$ such that
$Q(v_{k,\lambda})=0$ near $(1,1)$ with $k>1$. On the other hand,
under the assumptions of $V(x)$ and $W(x)$, it is easy to see that
$g(k,1)<0$ for any $k>1$. By the continuity, we can choose a pair
of $(k,\lambda)$ near $(1,1)$ with $k>1$ satisfies both
$Q(v_{k,\lambda})=0$ and $S_{\omega}(v_{k,\lambda})<0$. Letting
$u^*=v_{k,\lambda}$ for this $(k,\lambda)$, we get that
$S_{\omega}(u^*)<0$ and $Q(u^*)=0$.\hfill $\Box$

Proposition 5.1.2 means that $\mathcal{CM}$ is not empty and
$d_{\mathcal{M}}$ is well defined. Moreover, we have

{\bf Proposition 5.1.3.} \  {\it Assume that the conditions of
Theorem 4 hold. Then $d_{\mathcal{M}}>0$.}

{\bf Proof:} $u\in \Sigma\setminus\{0\}$ and $S_{\omega}(u)<0$
imply that
\begin{align}
&\qquad
2\omega\int_{\mathbb{R}^N}|u|^2dx+\int_{\mathbb{R}^N}[|\nabla
u|^2+V(x)|u|^2]dx\nonumber\\
&<\int_{\mathbb{R}^N}|u|^2f(x,|u|^2)dx+\int_{\mathbb{R}^N}(W\star|u|^2)|u|^2dx\nonumber\\
&\leq\frac{l+1}{l}\int_{\mathbb{R}^N}[c_1|u|^{2p_1+2}+c_2|u|^{2p_2+2}]dx\nonumber\\
&\quad+\|W_1\|_{L^{\infty}}\|u\|_{L^2}^4+C\|W_2\|_{L^q}\| \nabla
u\|_{L^2}^{\frac{N}{q}}\|u\|_{L^2}^{\frac{4q-N}{q}}.\label{97221}
\end{align}
Similar to (\ref{9814w1}) and (\ref{97171}), from (\ref{97221}), we
have
\begin{align}
2\omega\int_{\mathbb{R}^N}|u|^2dx+\int_{\mathbb{R}^N}[|\nabla
u|^2+V(x)|u|^2]dx\geq C>0.\label{97222}
\end{align}

On the other hand, if $Q(u)=0$, we have
\begin{align*}
&\qquad 2\int_{\mathbb{R}^N}|\nabla
u|^2dx-\int_{\mathbb{R}^N}(x\cdot\nabla V)
|u|^2dx\nonumber\\&=N\int_{\mathbb{R}^N}[|u|^2f(x,|u|^2)-F(x,|u|^2)]dx-\frac{1}{2}\int_{\mathbb{R}^N}
\{(x\cdot \nabla W)\star|u|^2\}|u|^2dx\nonumber\\
&\geq
Nl\int_{\mathbb{R}^N}F(x,|u|^2)dx+\frac{1}{2}\int_{\mathbb{R}^N}
\{(x\cdot \nabla W)\star|u|^2\}|u|^2dx,
\end{align*}
that is,
\begin{align}
&\qquad-\frac{1}{2}\int_{\mathbb{R}^N}F(x,|u|^2)dx+\frac{1}{4Nl}\int_{\mathbb{R}^N}
\{(x\cdot \nabla W)\star|u|^2\}|u|^2dx\nonumber\\
&\geq -\frac{1}{Nl}\int_{\mathbb{R}^N}|\nabla
u|^2dx+\frac{1}{2Nl}\int_{\mathbb{R}^N}(x\cdot\nabla V)
|u|^2dx.\label{9823w1}
\end{align}
Using (\ref{97251}), (\ref{9816x1}), (\ref{97252}), (\ref{97222})
and (\ref{9823w1}), we can get
\begin{align}
I_{\omega}(u)&=\omega\int_{\mathbb{R}^N}|u|^2dx+\frac{1}{2}\int_{\mathbb{R}^N}[|\nabla
u|^2+V(x)|u|^2-F(x,|u|^2)]dx-\frac{1}{4}\int_{\mathbb{R}^N}(W\star|u|^2)|u|^2dx\nonumber\\
&\geq
\omega\int_{\mathbb{R}^N}|u|^2dx+\frac{Nl-2}{2Nl}\int_{\mathbb{R}^N}|\nabla
u|^2dx+\frac{1}{2Nl}\int_{\mathbb{R}^N}[NlV(x)+(x\cdot \nabla V)]|u|^2dx\nonumber\\
&\qquad
-\frac{1}{4Nl}\int_{\mathbb{R}^N}\left\{[NlW+(x\cdot\nabla W)]\star|u|^2\right\}|u|^2dx\nonumber\\
&\geq
C\left(2\omega\int_{\mathbb{R}^N}|u|^2dx+\int_{\mathbb{R}^N}[|\nabla
u|^2+V(x)|u|^2]dx\right)\nonumber\\
&\geq C>0.\label{97223}
\end{align}
Consequently, $$ \qquad \qquad \qquad \qquad \qquad\qquad
d_{\mathcal{M}}=\inf_{\mathcal{CM}}I_{\omega}(u)>C>0. \qquad
\qquad \qquad \qquad \qquad  \qquad \qquad
 \Box$$

By the conclusions of Proposition 5.1.1 and Proposition 5.1.3,
we have
\begin{align}
d_{II}=\min\{d_{\mathcal{M}}, d_{\mathcal{N}}\}>0.\label{97224'}
\end{align}

Now we define the following manifolds:
\begin{align}
\mathcal{K}:&=\{u\in \Sigma\setminus\{0\}:I_{\omega}(u)<d_{II}, \
S_{\omega}(u)<0, \ Q(u)<0\},  \label{97225}\\
\mathcal{K}_+:&=\{u\in \Sigma\setminus\{0\}:I_{\omega}(u)<d_{II}, \
S_{\omega}(u)<0, \
Q(u)>0\},  \label{97226}\\
\mathcal{R}_+:&=\{u\in \Sigma\setminus\{0\}:I_{\omega}(u)<d_{II}, \
S_{\omega}(u)>0\}. \label{97227}
\end{align}
The following proposition will show some properties of
$\mathcal{K}$, $\mathcal{K}_+$ and $\mathcal{R}_+$:

{\bf Proposition 5.1.4} {\it Assume that the conditions of Theorem 4
hold. Then

(i) $\mathcal{K}$, $\mathcal{K}_+$ and $\mathcal{R}_+$ are not
empty.

(ii) $\mathcal{K}$, $\mathcal{K}_+$ and $\mathcal{R}_+$ are
invariant manifolds of (\ref{system0}).}

{\bf Proof:} (i) In order to prove  $\mathcal{K}$ is not empty, we
only need to find that there at least exists a $w\in \mathcal{K}$.
For $w^{\star}\in \Sigma\setminus \{0\}$ satisfies
$S_{\omega}(w^{\star})=0$ and $Q(w^{\star})=0$, letting
$w_{\rho}=\rho w^{\star}$ for $\rho
>0$, we have
\begin{align*}
S_{\omega}(w_{\rho})&=\rho^2
\int_{\mathbb{R}^N}\left\{2\omega|w^{\star}|^2+|\nabla
w^{\star}|^2+V(x)|
w^{\star}|^2\right\}dx\nonumber\\
&\quad-\int_{\mathbb{R}^N}\rho^2|w^{\star}|^2f(x,\rho^2|w^{\star}|^2)dx
-\rho^4\int_{\mathbb{R}^N}(W\star|w^{\star}|^2)|w^{\star}|^2dx,\\
Q(w_{\rho})&=\rho^2\int_{\mathbb{R}^N}\left(2|\nabla
w^{\star}|^2-(x\cdot \nabla V)|w^{\star}|^2\right)dx\nonumber\\
&\quad+N\int_{\mathbb{R}^N}[F(x,\rho^2|w^{\star}|^2)-\rho^2|w^{\star}|^2f(x,\rho^2|w^{\star}|^2)]dx\nonumber\\
&\quad +\frac{1}{2}\rho^4\int_{\mathbb{R}^N}\{(x\cdot \nabla W)\star|w^{\star}|^2\}|w^{\star}|^2dx,\\
 I_{\omega}(u_{\rho})&=\frac{1}{2}\rho^2
\int_{\mathbb{R}^N}\left\{2\omega|w^{\star}|^2+|\nabla
w^{\star}|^2+V(x)|
w^{\star}|^2\right\}dx\nonumber\\
&\quad -\frac{1}{2}\int_{\mathbb{R}^N}F(x,\rho^2|w^{\star}|^2)dx-
\frac{1}{4}\rho^4\int_{\mathbb{R}^N}(W\star|w^{\star}|^2)|w^{\star}|^2dx.
\end{align*}

Since $f(x,|w^*|^2)<f(x,\rho^2 |w^*|^2)$ and $\rho^2F(x,|w^*|^2)<F(x,\rho^2 |w^*|^2)$ for $\rho>1$
and from (\ref{9828w1}), we can obtain
\begin{align}
S_{\omega}(w_{\rho})&< \rho^2S_{\omega}(w^{\star})=0,\quad
Q(w_{\rho})< \rho^2Q(w^{\star})=0\label{97228}
\end{align}
for any $\rho>1$.  Noticing $d_{II}>0$, we also can choose
$\rho>1$ closing to 1 enough such that
\begin{align}
I_{\omega}(w_{\rho})<
\rho^2I_{\omega}(w^{\star})<d_{II}.\label{9823w3}
\end{align}
(\ref{97228}) and (\ref{9823w3}) means that $w_\rho \in
\mathcal{K}$. That is, $\mathcal{K}$ is not empty.

Similar to (\ref{97228}), we can obtain
\begin{align}
S_{\omega}(w_{\rho})> \rho^2S_{\omega}(w^{\star})=0.\label{97229}
\end{align}
for any $0<\rho<1$.  Noticing $d_{II}>0$, we also can choose
$0<\rho<1$ closing to 1 enough such that
$I_{\omega}(w_{\rho})<d_{II}$ by continuity, which implies that
$w_{\rho}\in \mathcal{R}_+$. That is, $\mathcal{R}_+$ is not
empty.

For $w^*\in \Sigma$ satisfies $S_{\omega}(w^*)<0$ and $Q(w^*)=0$,
letting $w_{\sigma}=\sigma w^*$ for $\sigma
>0$, we have
\begin{align*}
Q(w_{\sigma})&=\sigma^2\int_{\mathbb{R}^N}(2|\nabla
w^*|^2-(x\cdot \nabla V)|w^*|^2)dx\nonumber\\
&\quad
-\int_{\mathbb{R}^N}N[\sigma^2|w^*|^2f(x,\sigma^2|w^*|^2)-F(x,\sigma^2|w^*|^2)]dx\nonumber\\
&\quad+\frac{1}{2}\sigma^4\int_{\mathbb{R}^N}\left\{(x\cdot\nabla W)\star |w^*|^2\right\}|w^*|^2dx,\\
S_{\omega}(w_{\sigma})&=\sigma^2
\int_{\mathbb{R}^N}\left\{2\omega|w^*|^2+|\nabla w^*|^2+V(x)|
w^*|^2\right\}dx\nonumber\\
&\quad-\int_{\mathbb{R}^N}\sigma^2|w^*|^2f(x,\sigma^2|w^*|^2)dx
-\sigma^4\int_{\mathbb{R}^N} (W\star |w^*|^2)|w^*|^2dx,\\
I_{\omega}(w_{\sigma})&=\frac{1}{2}\sigma^2
\int_{\mathbb{R}^N}\left\{2\omega|w^*|^2+|\nabla w^*|^2+V(x)|
w^*|^2\right\}dx\nonumber\\
&\quad-\frac{1}{2}\int_{\mathbb{R}^N}F(x,\sigma^2|w^*|^2)dx-\frac{1}{4}\sigma^4\int_{\mathbb{R}^N}
(W\star |w^*|^2)|w^*|^2dx.
\end{align*}

Since $\phi(\sigma)=Q(w_{\sigma})$ is a smooth function  of $\sigma$
and $Q(w^*)=0$, we have  $\phi(1)=0$.  If $\phi'(1)\neq 0$,
then there exists a $\sigma_0>0$ such that
$Q(u_{\sigma})=\phi(\sigma)>0$ for $\sigma\in (1,\sigma_0)$ if
$\sigma_0>1$(or $\sigma\in (\sigma_0,1)$ if $\sigma_0<1$). By
continuity, we can choose such $\sigma_0$ closing to 1 enough such
that $S_{\omega}(w_{\sigma})<0$ and $I_{\omega}(w_{\sigma})<d_{II}$
for $\sigma\in (1,\sigma_0)$ if $\sigma_0>1$(or $\sigma\in
(\sigma_0,1)$ if $\sigma_0<1$). That is, $w_{\sigma}\in
\mathcal{K}_+$  and $\mathcal{K}_+$ is not empty.

If $\phi'(1)=0$, from $\phi(1)=0$ and $\phi'(1)=0$, we can respectively obtain
\begin{align*}
&\qquad -N\int_{\mathbb{R}^N}[|w^*|^2f(x,|w^*|^2)-F(x,|w^*|^2)]dx\nonumber\\
&=-N\int_{\mathbb{R}^N}|w^*|^4f'_s(x,|w^*|^2)dx+\frac{1}{2}\int_{\mathbb{R}^N}\left\{(x\cdot\nabla
W)\star |w^*|^2\right\}|w^*|^2dx
\end{align*}
and
\begin{align*}
Q(w^*)&=\int_{\mathbb{R}^N}\left(2|\nabla w^*|^2-(x\cdot \nabla
V)|w^*|^2
-N|w^*|^4f'_s(x,|w^*|^2)\right)dx\nonumber\\
&\quad+\int_{\mathbb{R}^N}\left\{(x\cdot\nabla W)\star
|w^*|^2\right\}|w^*|^2dx.
\end{align*}

Letting $w_{\sigma}=\sigma w^*$, we have
\begin{align}
Q(w_{\sigma})&=\sigma^2\int_{\mathbb{R}^N}\left(2|\nabla
w^*|^2-(x\cdot \nabla V)|w^*|^2
-N|w^*|^4f'_s(x,\sigma^2|w^*|^2)\right)dx\nonumber\\
&\quad+\sigma^4\int_{\mathbb{R}^N}\left\{(x\cdot\nabla W)\star
|w^*|^2\right\}|w^*|^2dx\nonumber\\
&>\sigma^2\int_{\mathbb{R}^N}\left(2|\nabla w^*|^2-(x\cdot \nabla
V)|w^*|^2
-N|w^*|^4f'_s(x,|w^*|^2)\right)dx\nonumber\\
&\quad+\sigma^4\int_{\mathbb{R}^N}\left\{(x\cdot\nabla W)\star
|w^*|^2\right\}|w^*|^2dx\nonumber\\
&=\sigma^2
Q(w^*)+(\sigma^4-\sigma^2)\int_{\mathbb{R}^N}\left\{(x\cdot\nabla
W)\star |w^*|^2\right\}|w^*|^2dx>0\label{9722w1}
\end{align}
for $0<\sigma<1$. By continuity, we can choose such $\sigma$
closing to 1 enough such that $S_{\omega}(w_{\sigma})<0$ and
$I_{\omega}(w_{\sigma})<d_{II}$. That is to say, $w_{\sigma}\in
\mathcal{K}_+$  and $\mathcal{K}_+$ is not empty.

(ii) In order to prove that $\mathcal{K}$ is the invariant
manifold of (\ref{system0}),  we need to show that: If $u_0\in
\mathcal{K}$, then solution $u(x,t)$ of (\ref{system0}) satisfies
$u(x,t)\in \mathcal{K}$ for any $t\in [0,T)$.

Assume that $u(x,t)$ is a solution of (\ref{system0}) with $u_0\in
\mathcal{K}$. Then we can obtain
\begin{align}
I_{\omega}(u(\cdot,t)&=E(u(\cdot,t))+\omega
\|u(\cdot,t)\|_2^2=E(u_0)+\omega\|u_0\|^2_2=I_{\omega}(u_0)<d_{II}\label{881018}
\end{align}
for $t\in [0,T)$. Next we prove that $S_{\omega}(u(\cdot,t))<0$ for $t\in [0,T)$.
Otherwise, by continuity, there exists a $t_0\in (0,T)$ such that
$S_{\omega}(u(\cdot,t_0))=0$ because of $S_{\omega}(u_0)<0$. Since
$\|u(\cdot,t)\|_2^2=\|u_0\|_2^2$ and $u_0\in \Sigma\setminus\{0\}$,
it is easy to see that $u(\cdot,t_0)\in \Sigma\setminus\{0\}$. By
the definitions of $d_{\mathcal{N}}$ and $d_{II}$, we know that
$I_{\omega}(u(\cdot,t_0))\geq d_{\mathcal{N}}\geq d_{II}$, which is
a contradiction to $I_{\omega}(u(\cdot,t))<d_{II}$ for $t\in [0,T)$.
Hence $S_{\omega}(u(\cdot,t))<0$ for all $t\in [0,T)$.

Now we only need to prove that $Q(u(\cdot,t))<0$ for $t\in [0,T)$.
Otherwise, since $Q(u_0)<0$, there exists a $t_1\in (0,T)$ such
that $Q(u(\cdot,t_1))=0$ by continuity. And
$S_{\omega}(u(\cdot,t_1))<0$ means that $u(\cdot,t_1)\in
\mathcal{CM}$. By the definitions of $d_{\mathcal{M}}$ and
$d_{II}$, we obtain $I_{\omega}(u(\cdot,t_1))\geq
d_{\mathcal{M}}\geq d_{II}$, which is a contradiction to
$I_{\omega}(u(\cdot,t))<d_{II}$ for $t\in [0,T)$. Hence
$Q(u(\cdot,t))<0$ for all $t\in [0,T)$.

By the discussions above, we know that:  $u(x,t)\in
\mathcal{K}$ for any $t\in [0,T)$ if $u_0\in \mathcal{K}$, which
means that $\mathcal{K}$ is the invariant manifold of
(\ref{system0}).

Similarly, we can prove that $\mathcal{K}_+$ and $\mathcal{R}_+$ are
also invariant manifolds of (\ref{system0}). \hfill $\Box$

{\bf Remark 5.1.1.} By the definitions of $d_{II}, d_{\mathcal{N}},
d_{\mathcal{M}}$, $\mathcal{K}$, $\mathcal{K}_+$ and
$\mathcal{R}_+$, it is easy to see that
$$
\qquad \qquad \qquad \{u\in \Sigma\setminus\{0\}:\
I_{\omega}(u)<d_{II}\}=\mathcal{K}\cup
\mathcal{K}_+\cup\mathcal{R}_+.
$$

\subsection{The proof of Theorem 4}
\qquad The proof of Theorem 4 depends on the following two lemmas.

{\bf Lemma 5.2.1.}  {\it Assume that the conditions of Theorem 4
hold. Then the solutions of (\ref{system0}) with $u_0\in
\mathcal{K}$ will blow up in finite time.}

{\bf Proof:} Since $u_0\in \mathcal{K}$ and $\mathcal{K}$ is the
invariant manifold of (\ref{system0}), we have $Q(u(x,t))<0$,
$S_{\omega}(u(x,t))<0$ and $I_{\omega}(u(x,t))<d_{II}$.

Under the conditions of Theorem 4, we have $J''(t)=4Q(u)<0$ and
$J'(0)<0$. By the results of Proposition 2.2, the solution
$u(x,t)$ will blow up in finite time. The conclusion of this lemma
is true.\hfill $\Box$

On the other hand, we have a parallel result on global existence.

 {\bf Lemma 5.2.2.} {\it Assume that the conditions of Theorem 4
hold. If $u_0\in \mathcal{K}_+$ or $u_0\in\mathcal{R}_+$, then the
solutions of (\ref{system0}) exists globally.}

{\bf Proof:} Case 1: Assume that $u(x,t)$ is a solution of
(\ref{system0}) with $u_0\in \mathcal{K}_+$. Since $\mathcal{K}_+$
is a invariant manifold of (\ref{system0}), we know that
$u(\cdot,t)\in \mathcal{K}_+$, which means that
$I_{\omega}(u(\cdot,t))<d_{II}$ and $Q(u(\cdot,t))>0$.
$Q(u(\cdot,t))>0$ and (\ref{97251}) imply that
\begin{align}
&\quad 2\int_{\mathbb{R}^N}|\nabla
u|^2dx-\int_{\mathbb{R}^N}(x\cdot \nabla V)|u|^2dx
\nonumber\\
&\geq Nl\int_{\mathbb{R}^N}
F(x,|u|^2)dx-\frac{1}{2}\int_{\mathbb{R}^N}\left\{(x\cdot\nabla
W)\star |u|^2\right\}|u|^2dx.\label{97231}
\end{align}
By the definition of $I_{\omega}(u)$ and using (\ref{97231}), we
have
\begin{align}
d_{II}>I_{\omega}(u(\cdot,t))&=\omega\int_{\mathbb{R}^N}|u|^2dx
+\frac{1}{2}\int_{\mathbb{R}^N}[|\nabla u|^2+V(x)|u|^2]dx\nonumber\\
&\quad
-\frac{1}{2}\int_{\mathbb{R}^N}F(x,|u|^2)dx-\frac{1}{4}\int_{\mathbb{R}^N}
(W\star |u|^2)|u|^2dx\nonumber\\
&\geq
\omega\int_{\mathbb{R}^N}|u|^2dx+\frac{Nl-2}{2Nl}\int_{\mathbb{R}^N}|\nabla
u|^2dx\nonumber\\
&\quad+\int_{\mathbb{R}^N}\frac{NlV(x)+(x\cdot \nabla V)}{2Nl}|u|^2dx\nonumber\\
&\quad
-\frac{1}{4Nl}\int_{\mathbb{R}^N}\left\{[NlW+(x\cdot\nabla W)]\star|u|^2\right\}|u|^2dx\nonumber\\
 &\geq
C\left(\int_{\mathbb{R}^N}|u|^2dx+\int_{\mathbb{R}^N}|\nabla
u|^2dx+\int_{\mathbb{R}^N}V(x)|u|^2dx\right).\label{97232}
\end{align}
(\ref{97232}) means that $u(x,t)$ exists globally.

Case 2:  Assume that $u(x,t)$ is a solution of (\ref{system0})
with $u_0\in\mathcal{R}_+$. Since $\mathcal{R}_+$ is also a
invariant manifold of (\ref{system0}), we know that $u(x,t), \in
\mathcal{R}_+$, which means that $I_{\omega}(u(\cdot,t))<d_{II}$
and $S_{\omega}(u(\cdot,t))>0$. Since $S_{\omega}(u)>0$, we can
get
\begin{align}
&\quad\omega \|u\|_2^2+\frac{1}{2}\int_{\mathbb{R}^N}(|\nabla
u|^2+V(x)|u|^2)dx\nonumber\\
&>\frac{1}{2}\int_{\mathbb{R}^N}f(x,|u|^2)|u|^2dx+\frac{1}{2}\int_{\mathbb{R}^N}(W\star|u|^2)|u|^2dx\nonumber\\
&\geq \min(l+1,2)\left(\frac{1}{2}\int_{\mathbb{R}^N}F(x,|u|^2)dx
+\frac{1}{4}\int_{\mathbb{R}^N}(W\star|u|^2)|u|^2dx\right).\label{9823wx1}
\end{align}
From (\ref{9823wx1}), we can obtain
\begin{align}
I_{\omega}(u)&=\omega\|u\|_2^2+\frac{1}{2}\int_{\mathbb{R}^N}[|\nabla
u|^2+V(x)|u|^2-F(x,|u|^2)]dx-G(|u|^2)\nonumber\\
&\geq
\min\left(\frac{l}{(l+1)},\frac{1}{2}\right)\left(\omega\|u\|_2^2+\frac{1}{2}\int_{\mathbb{R}^N}[|\nabla
u|^2+V(x)|u|^2]dx\right).\label{9823wx3}
\end{align}
(\ref{9823wx3}) implies that the solution $u(x,t)$ exists
globally.\hfill $\Box$

{\bf The proof of Theorem 4:} By the results of Lemma 5.2.1, Lemma
5.2.2, we know that Theorem 4 is right.\hfill $\Box$

As a corollary of Theorem 4, we obtain a sharp threshold for the blowup in finite time and
global existence of the solution of
(\ref{9829w1}) as follows

{\bf Corollary 5.1.} {\it Assume that $f(x,|u|^2)\equiv 0$, $V(x)\equiv 0$, $W(x)>0$ for all $x\in \mathbb{R}^N$,
$W$ is even and $W\in L^{\infty}(\mathbb{R}^N)+L^q(\mathbb{R}^N)$
with some $q>\frac{N}{4}$. Suppose further that there exists $l$
satisfying $2<Nl$ and
\begin{align*}
 NlW(x)+(x\cdot \nabla W)\leq 0.
\end{align*}  If $u_0\in H^1(\mathbb{R}^N)$, $|x|u_0\in L^2(\mathbb{R}^N)$ and $I_{\omega}(u_0)=\omega
\|u_0\|_2^2+E(u_0)<d_{II}$, then the solution of (\ref{9829w1})
blows up in finite time if and only if $u_0\in \mathcal{K}$.}

{\bf Remark 5.2.1.} A typical example of (\ref{9829w1}) is
\begin{align} \label{Hartree} \left\{\begin{array}{ll}&-iu_t=\Delta u+(|x|^{-K}\star |u|^2)u, \quad x\in\mathbb{R}^N, \quad t>0,\\
&u(x,0)=u_0(x), \quad x\in\mathbb{R}^N,
\end{array}
\right. \end{align}
which is also a special case of (\ref{system0}) with $V(x)\equiv 0$,
$f(x,|u|^2)\equiv 0$ and $W(x)=|x|^{-K}$ with $2<Nl<K<\frac{N}{q}<4$.
Letting $W=W_1+W_2$ with
$$
W_1(x)=\{^{0,\qquad \ |x|\leq 1,}_{|x|^{-K},\quad |x|>1}\quad {\rm
and}\quad  W_2(x)=\{^{|x|^{-K},\quad |x|\leq 1,}_{0,\qquad \
|x|>1,}
$$
we can see that $W_1\in L^{\infty}(\mathbb{R}^N)$ and $W_2\in
L^q(\mathbb{R}^N)$ with some $\frac{N}{4}<q<\frac{N}{2}$.
Corollary 5.1 gives the sharp threshold for blowup and global existence of the solution to (\ref{Hartree}).

We will give some examples of $V(x)$, $f(x,|u|^2)$ and $W(x)$. It is easy to verify that they satisfy the conditions of Theorem 4.

Example 1. $V(x)=|x|^2$, $W(x)=a|x|^{-K}$ with
$2<Nl<K<\frac{N}{q}<4$ for $x\in \mathbb{R}^N$ and
$f(x,|u|^2)=b|u|^{2p_1}+c|u|^{2p_2}$ with $a\geq 0$, $b>0$,
$c>0$ and $p_2>p_1>\frac{2}{N}$.

Example 2. $V(x)=|x|^2$, $W(x)=a|x|^{-K}$ with
$2<Nl<K<\frac{N}{q}<4$ for $x\in \mathbb{R}^N$ and
$f(x,|u|^2)=c|u|^{2q_1}+d|u|^{2q_2}$ with $a\geq 0$, $c$ is a
real number, $d>0$ and $q_2>\frac{2}{N}$, $q_2>q_1>0$.

Example 3. $V(x)=\frac{|x|^2}{1+|x|^2}$, $W(x)=a|x|^{-K}$ with
 $2<Nl<K<\frac{N}{q}<4$ for $x\in \mathbb{R}^N$ and $f(x,|u|^2)=b|u|^{2p}\ln(1+|u|^2)$
 with $a\geq 0$, $b>0$ and $p>\frac{2}{N}$.

\end{document}